\documentclass[10pt, final, journal, letterpaper, twoside, twocolumn]{IEEEtran}
\usepackage{graphicx}
\usepackage{subcaption}
\usepackage{epsfig} 
\usepackage{mathptmx} 
\usepackage{times} 
\usepackage{amsmath} 
\usepackage{amssymb}
\usepackage{enumerate}
\usepackage{mathrsfs}
\usepackage{euscript}
\usepackage{psfrag}
\usepackage[section]{placeins} 
\usepackage{amsbsy}
\usepackage{textcomp}
\usepackage{blkarray}
\usepackage{mathdots}
\usepackage[usenames,dvipsnames]{color}
\usepackage{mdframed}
\usepackage{caption}
\usepackage{subcaption}
\usepackage[T1]{fontenc}

\newtheorem{theorem}{Theorem}[section]

\newtheorem{lemma}[theorem]{Lemma}
\newtheorem{example}[theorem]{Example}
\newtheorem{remark}[theorem]{Remark}

\newtheorem{problem}[theorem]{Problem}

\begin{document}
 \title{Steering nonholonomic integrator using orthogonal polynomials}
 \author{Pragada Shivaramakrishna, A. Sanand Amita Dilip \thanks{Pragada Shivaramakrishna is with the Department of Aerospace Engineering, 
 Indian Institute of Technology Kharagpur, India. {\tt\small shivaramkratos@gmail.com, shivaram@iitkgp.ac.in}}\thanks{
  A. Sanand Amita Dilip is with the Department of Electrical Engineering, Indian Institute of Technology Kharagpur, India. 
 {\tt\small sanand@ee.iitkgp.ac.in}} 
}
 \maketitle
 \begin{abstract}
  We consider minimum energy optimal control problem with time dependent Lagrangian on the nonholonomic integrator and and find the analytical solution 
  using Sturm-Liouville theory. Furthermore, we also consider the minimum energy problem on the Lie group $\mathbb{SO}(3)$ with time dependent Lagrangian. 
  We show that the steering of nonholonomic integrator and generalized nonholonomic integrator 
  can be achieved by using various families of orthogonal polynomials such as Chebyshev, Legendre 
  and Jacobi polynomials apart from trigonometric polynomials considered in the literature. 
  Finally, we show how to  
  find sub-optimal inputs using elements from a family of orthogonal functions when the cost function is given by the $\mathcal{L}_1$ norm of the input.
 \end{abstract}
 \section{Introduction}
 
 In this paper, we discuss the minimum energy steering problem on the nonholonomic integrator model where the cost function is given by time dependent Lagrangians. 
 The proposed approach reveals a connection between optimal control, Sturm-Liouville theory and orthogonal polynomials which can be used to study 
 optimal control problems with time dependent Lagrangians for general models with higher order nonlinear effects.  
 
Motion planning is an important problem in control theory and its applications such as robotics, spacecraft attitude and guidance control, vehicle control 
and so on (\cite{murray},\cite{murrayli},\cite{henninger},\cite{henninger1},\cite{biggs},\cite{paden},\cite{sun}). Algorithms 
 for motion planning has always been an active area of research (\cite{murray},\cite{sastry},\cite{walsh1},\cite{biggs},\cite{henninger} and the references 
 therein). Many mechanical systems arising in engineering 
 can be modeled as nonholonomic systems e.g., vehicle models, robot arm manipulator, spacecraft models and so on 
 (\cite{murray},\cite{murrayli},\cite{sastry},\cite{biggs},\cite{Spindler},\cite{Spindler1},\cite{Spindler2}). Efficient 
 motion planning becomes necessary in the presence of obstacles and uncertainties in path planning and tracking problems \cite{sun}. 
 The nonholonomic integrator provides the simplest model of nonholonomic systems and serves as a  
 prototype to understand general nonlinear systems \cite{broc2}. This model can be generalized to model many mechanical systems arising in engineering 
 \cite{murray},\cite{sastry},\cite{broc2},\cite{broc1}. A slightly more general version of the nonholonomic integrator model involves chained form 
 systems \cite{murray},\cite{sastry}. Cooperative control of nonholonomic mobile agents where each agent has has chained form dynamics is studied 
 in \cite{dong}. 

 The nonholonomic integrator is defined by  
\begin{eqnarray}
  \dot{x}_1=u_1\;, \dot{x}_2=u_2\;,\dot{x}_3=x_1u_2-x_2u_1.\label{noholintg}
 \end{eqnarray}
 We are interested in steering this system from a given initial state to the desired final state in finite time using optimal/sub-optimal inputs with respect 
 to the given cost function. 
 The nonholonomic integrator model or the Brockett integrator has an extremely rich structure as can be seen by the literature on this topic 
 (\cite{sastry},\cite{broc2},\cite{broc1} and the references therein). We focus on the steering aspects only and do not consider the stabilization problem in here. 
 Amongst its 
 various properties, Murray and Sastry (\cite{murray}) explored the steering of the nonholonomic integrator using sinusoidal input of appropriate frequencies 
 to give sub-optimal inputs. These methods were then extended for more general nonholonomic integrator models as well such as the generalized 
 nonholonomic integrator and the chained form systems (\cite{murray}). The orthogonality property of 
 sinusoids was exploited in this approach. It was shown in \cite{lib1} that the following unicycle system 
 \begin{eqnarray}
  \dot{x}_1=u_1\cos\theta\;, \dot{x}_2=u_2\sin \theta\;,\dot{\theta}=u_2.\nonumber
 \end{eqnarray}
 can be converted into the form $(\ref{noholintg})$. Thus, the steering of unicycle models is reduced to the steering of the nonholonomic integrator. 
 
 Another very commonly used model for control of rigid bodies such as spacecraft is that of a nonholonomic control system on the Lie group $\mathbb{SO}(3)$. 
  Optimal control problems on Lie groups have been studied in \cite{henninger},\cite{henninger1},\cite{biggs},\cite{walsh1},\cite{walsh},\cite{Spindler},\cite{Spindler1},\cite{Spindler2} 
 and the references therein. We refer the reader to \cite{Spindler} for a tutorial type treatment of optimal control on Lie groups.   
 Optimal control problems on the Lie group $\mathbb{SO}(3)$ can be identified with optimal control problems involving 
 spacecraft attitude control and so on. Commonly studied problems involve underactuated spacecraft attitude control \cite{henninger1},\cite{Spindler1}.
 
 The minimum energy optimal control problems considered in the above mentioned references assume 
 time independent Lagrangian. We consider more general time dependent Lagrangians and solve the optimal control problem analytically 
 by exploiting orthogonality properties of families of orthogonal functions such 
 as Chebyshev polynomials and Sturm-Liouville theory. This can also be considered as an inverse optimal control problem where one finds an appropriate 
 cost function when the input is given by orthogonal polynomials. 
 We then apply these ideas to an optimal control problem on the Lie group $\mathbb{SO}(3)$ with time dependent Lagrangian. 
  We also extend the steering 
 algorithm of \cite{murray} using trigonometric functions to more general families of orthogonal polynomials in the case of the nonholonomic integrator. \\
 {\bf Organization}: In the next section, we give some preliminaries on nonlinear controllability, optimal control using sinusoids, Sturm-Liouville equations 
 and orthogonal polynomials. This is followed by a section where we analytically solve minimum energy optimal control problem on the nonholonomic integrator 
 with time dependent Lagrangian. 
 In Section 
 $IV$, we extend the steering algorithm of \cite{murray} to families of orthogonal polynomials to give sub-optimal inputs. We then extend our approach 
 to under-actuated control of Lie group $\mathbb{SO}(3)$. Then, we show how to obtain sub-optimal control input for different families of orthogonal 
 polynomials when the cost function on the control inputs is given by the $\mathcal{L}_1-$norm. \\
 {\bf Notation}: Vectors are denoted by boldface letters, matrices are denoted by capital letters and scalars  
 by small-face letters. A vector valued function 
 $\bold{f}=[f_1,\ldots,f_n]$ of variables $x_1,\ldots,x_n$ is denoted by $\bold{f}(\bold{x})$. Its Jacobian matrix is denoted by $\frac{\partial \bold{f}}{\partial \bold{x}}$ whose $(i,j)-$th 
 entry is $\frac{\partial f_i}{\partial x_j}$. The time derivative of a function or a variable say $x_i$ is sometimes denoted by $\dot{x}_i$. Elements of Lie 
 groups are denoted by small-case letters although they are be represented by matrices. The special orthogonal group of $n\times n$ matrices (the subset of 
 orthogonal matrices with the determinant equal to one) 
 is denoted by $\mathbb{SO}(n)$. The two norm or the $\mathcal{L}_2$ norm on function spaces is 
 is represented by $\|.\|$, the one norm or $\mathcal{L}_1$ norm by $|.|$ and general $\mathcal{L}_p$ norm by $\|.\|_p$. 
 Lie bracket operation is denoted by $[.,.]$ and ad$(U)(E):=[U,E]=UE-EU$ where $U,E$ are $n\times n$ matrices is the matrix commutator. The symbol $\circ$ 
 denotes the composition of functions or the action of an element on elements of some other set for example, the action of a dual vector from dual vector 
 space on the elements of the primal vector space. 
 
 Chebyshev polynomials of the first and the second kind are denoted by $T_n(t)$ and $U_n(t)$ respectively where $n\in \mathbb{N}\cup 0$. Jacobi 
 polynomials are denoted by $P_n^{(\alpha,\beta)}(t)$ and Legendre polynomials by $P_n(t)$. The symbol $\delta_{mn}$ denotes the Kronecker-delta function. 
 \section{Preliminaries}
  {\bf Basics on controllability and Lie brackets}: 
  The formal definitions of manifolds, tangent spaces, vector fields, Lie groups and Lie algebras can be found for example in \cite{sastry},\cite{bullo}. 
 The Lie bracket between two vector fields $\bold{f},\bold{g}$ is defined as 
 \begin{equation}
  [\bold{f},\bold{g}]=\frac{\partial \bold{g}}{\partial \bold{x}}\bold{f}-\frac{\partial \bold{f}}{\partial \bold{x}}\bold{g}.\nonumber
 \end{equation}
The Lie algebra generated by the set of vector fields $\bold{f},\bold{g}_1,\ldots,\bold{g}_m$ consists of a vector space spanned by 
$\bold{f},\bold{g}_1,\ldots,\bold{g}_m$ which is closed under Lie bracket operation.

 A general nonlinear control system $\dot{\bold{x}}=\bold{f}(\bold{x},\bold{u})$ is said to be controllable if given an arbitrary initial condition $\bold{x}_0$, 
 one can drive the state $\bold{x}(t)$ to any arbitrary $\bold{x}_f$ in finite time using some input $\bold{u}$. Consider affine nonlinear 
 control system of the form $\dot{\bold{x}}=\bold{f}(\bold{x})+G(\bold{x})\bold{u}$ where columns of the matrix $G=[\bold{g}_1,\ldots,\bold{g}_m]$ form 
 smooth, linearly independent vector 
 fields. Checking controllability locally for such systems is equivalent to checking the Lie algebraic rank condition of the Lie algebra formed by the vector fields 
 $\bold{f}$ and $\bold{g}_i$, $i=1,\ldots,m$. We refer the reader to  \cite{sastry},\cite{bullo},\cite{isi},\cite{nij},\cite{suss} for the formal definitions and 
 different versions of controllability such as local/global controllability, accessibility, small time local controllability and so on. It turns out 
 that the nonholonomic integrator is globally controllable. Affine nonlinear control systems $\dot{\bold{x}}=\bold{f}(\bold{x})+G(\bold{x})\bold{u}$ 
 where $\bold{f}=\bold0$ are called drift-free or nonholonomic systems. 
 
 Consider an optimal control problem 
 \begin{eqnarray}
  &\mbox{minimize}_{\bold{u}}& J=\int_{t_0}^{t}L(\bold{x},\bold{u})dt\nonumber\\
  &\mbox{subject to}& \dot{\bold{x}}=\bold{f}(\bold{x},\bold{u}).\label{optctrl}
 \end{eqnarray}
The term $L(\bold{x},\bold{u})$ is called Lagrangian. One obtains an augmented Lagrangian $L_a(\bold{x},\bold{u},\bold{p}):=
L(\bold{x},\bold{u})+\bold{p}^{\mathsf{T}}(\dot{\bold{x}}-\bold{f}(\bold{x},\bold{u}))$ where entries of $\bold{p}$ are called Lagrange multipliers. Then, the Euler-Lagrange equation which is the necessary 
first order optimality condition for $(\ref{optctrl})$ is given by 
\begin{eqnarray}
 \frac{d}{dt}\Bigg(\frac{\partial L_a}{\partial \dot{\bold{x}}}\Bigg)=\frac{\partial L_a}{\partial \bold{x}}.\label{eullag}
\end{eqnarray}
Let $H(\bold{x,\bold{u},\bold{p}}):=\bold{p}^{\mathsf{T}}\bold{f}-L(\bold{x},\bold{u})$ be the Hamiltonian. Then, the first order optimality 
conditions using $H$ are given by the following Hamilton's equations 
\begin{eqnarray}
 \dot{\bold{x}}=\frac{\partial H}{\partial \bold{p}}, \dot{\bold{p}}=-\frac{\partial H}{\partial \bold{x}}, \frac{\partial H}{\partial \bold{u}}=0.\label{hamilteq}
\end{eqnarray}
In addition, there are boundary conditions depending on the given control problem. We refer the reader to \cite{lib} for further details on
variational approach to optimal control. For a background on calculus of variations, we refer the reader to \cite{gelf}.  
\\
  {\bf Optimal inputs and steering using sinusoids}: 
Consider the following example from \cite{sastry}. 
 \begin{example}[\cite{sastry}]\label{sasex}
%
 For system $(\ref{noholintg})$, we want to find the minimum energy input to drive the state from the origin to a specified point $(0,0,a)$ from $t=0$ to $t=1$. The cost function 
 is $J=\int_0^1(u_1^2+u_2^2)dt$ subject to system dynamics. 
 Using system equations to eliminate $u_1$ and $u_2$, we obtain the cost function $\int_0^1(\dot{x}_1^2+\dot{x}_2^2)dt$ subject to 
 $\dot{x}_3-x_1\dot{x}_2+x_2\dot{x}_1=0$. Therefore, 
 the augmented cost function is 
 \begin{eqnarray}
  J_a=\int_0^1(\dot{x}_1^2+\dot{x}_2^2+p(t)(\dot{x}_3-x_1\dot{x}_2+x_2\dot{x}_1))dt.\nonumber
 \end{eqnarray}
Applying the first order necessary conditions from calculus of variations (Euler-Lagrange equation $(\ref{eullag})$), we obtain 
\begin{eqnarray}
 -p(t)\dot{x}_2&=&\frac{d}{dt}(2\dot{x}_1+p(t)x_2)\nonumber\\
 p(t)\dot{x}_1&=&\frac{d}{dt}(2\dot{x}_2-p(t)x_1)\nonumber\\
 0&=&\frac{d}{dt}p(t).\nonumber
\end{eqnarray}
Thus, $p(t)=c$ and we obtain the following second order system of equations 
\begin{eqnarray}
 \ddot{x}_1+c\dot{x}_2&=&0\nonumber\\
 \ddot{x}_2-c\dot{x}_1&=&0.\nonumber
\end{eqnarray}
Now using $\dot{x}_1=u_1$ and $\dot{x}_2=u_2$, we have the following first order ode
\begin{eqnarray}
 \dot{\left[ \begin{array}{c} {u}_1\\u_2\end{array} \right]}&=&\left[ \begin{array}{cc} 0&-c\\c&0\end{array} \right]
 \left[ \begin{array}{c} {u}_1\\u_2\end{array} \right]\nonumber\\
 \Rightarrow \left[ \begin{array}{c} {u}_1\\u_2\end{array} \right]&=&\left[ \begin{array}{cc} \cos ct&-\sin ct\\\sin ct& \cos ct\end{array} \right]
 \left[ \begin{array}{c} {u}_1(0)\\u_2(0)\end{array} \right].\nonumber
\end{eqnarray}
We need to find $\bold{u}(0)$ and $c$ using initial and final conditions. Let's write $\dot{\bold{u}}=H\bold{u}$ for first order equations in $u_1,u_2$. 
Hence, $\bold{u}(t)=e^{Ht}\bold{u}(0)$. Note that $e^{Ht}$ is orthogonal, hence, the norm of $\|\bold{u}(t)\|=\|\bold{u}(0)\|$ remains constant for 
all time. 
 From the terminal conditions, it follows that $c=2n\pi$ where $n=0,\pm 1,\pm 2,\ldots$. 
Suppose $a>0$, then the cost is 
minimum when $n=1$ and $\|\bold{u}\|=2\pi a$ 
with the direction of $\bold{u}$ being arbitrary.

For an arbitrary terminal time $T$, 
it turns out that $cT=2n \pi$. 
Thus, for $n=1$, $c=\frac{2\pi}{T}$ and  
\begin{equation}
 \left[ \begin{array}{c} {u}_1(t)\\u_2(t)\end{array} \right]=\left[ \begin{array}{cc} \cos \frac{2\pi}{T} t&-\sin \frac{2\pi}{T} t\\\sin \frac{2\pi}{T} t& 
 \cos \frac{2\pi}{T} t\end{array} \right]
 \left[ \begin{array}{c} {u}_1(0)\\u_2(0)\end{array} \right].\nonumber
\end{equation}
Let $u_i(0)=\sqrt{\frac{ca}{2}}$, $i=1,2$. 
Therefore, with sinusoidal inputs of appropriate frequencies, one can always 
steer the system from the origin to any point $(0,0,a)$ in time $T$. The frequencies are chosen depending upon the terminal time $T$ so that for $x_1$ and 
$x_2$, we are integrate the sinusoids over the full period. 

 \end{example}
 
 \begin{remark}
  Notice that the orthogonality of trigonometric polynomials is crucial when we want to steer the system along the $x_3-$direction. This leads to a question 
whether one can use other families of orthogonal functions such as Legendre polynomials, Chebyshev polynomials and so on 
for steering along the $x_3-$direction and whether such orthogonal functions serve as optimal inputs for 
an appropriate cost function. 
 \end{remark}
  {\bf Sturm-Liouville equations and orthogonal polynomials}: 
  We briefly mention about Sturm-Liouville equations and orthogonal polynomials from \cite{arf}.
  Let $P=P(t)>0$ and $Q=Q(t)$ be two given functions such that $Q$ is continuous and $P$ is continuously differentiable. Then the following ode 
  \begin{equation}
   \frac{d}{dt}\Bigg(P\frac{dy}{dt}\Bigg)+Qy=-\lambda y\label{sleqn}
  \end{equation}
is called the Sturm-Liouville equation (\cite{arf}). Notice that it is an eigenvalue problem where $\lambda$ is an eigenvalue and $y$ is the corresponding 
eigenfunction. The Sturm-Liouville equation can also be represented in operator form as 
\begin{equation}
 \mathcal{L}y=\frac{d}{dt}\Bigg(P\frac{dy}{dt}\Bigg)+Qy.\label{slop}
\end{equation}
The Sturm-Liouville operator can also be defined as  
\begin{equation}
 \mathcal{L}y=\frac{1}{w(t)}\frac{d}{dt}\Bigg(P\frac{dy}{dt}\Bigg)+Qy\label{slop1}
\end{equation}
where $w$ is some weight factor. 
The eigenvalue equation is $\mathcal{L}y=-\lambda y$. 
If $f,g$ are solutions of $(\ref{slop1})$, then 
\begin{equation}
 \int_{a}^{b}f(t)g(t)w(t)dt=0\nonumber
\end{equation}
which implies that $f,g$ are orthogonal polynomials with respect to the weight factor $w$. 

The solutions of the following Sturm-Liouville differential equation 
\begin{equation}
 \frac{d}{dt}\Bigg((1-t^2)\frac{dy}{dt}\Bigg)+n(n+1)y=0.\label{Legendre}
\end{equation}
 with $n=0,1,2,\ldots$ are called Legendre polynomials. Let $y_n(t)=P_n(t)$ denote the family of solutions for $n=0,1,2,\ldots$. It turns out that 
 $\int_1^1P_m(t)P_n(t)=\frac{2}{2n+1}\delta_{mn}$ i.e., these polynomials are orthogonal in the interval $[-1,1]$.

 The Sturm-Liouville equation for Chebyshev polynomials of the first kind are given by choosing $p(t)=\sqrt{1-t^2}$, $w(t)=\frac{1}{\sqrt{1-t^2}}$ and 
 $\lambda=n^2$ (\cite{arf}) which are denoted by $T_n(t)$. Chebyshev polynomials of the second kind are given by choosing $p(t)=(1-t^2)^{\frac{3}{2}}$, $w(t)=\sqrt{1-t^2}$ and 
 $\lambda=n(n+2)$ (\cite{arf}) and are denoted by $U_n(t)$.
 The Sturm-Liouville differential equation for Jacobi polynomials is given by
\begin{eqnarray}
 \frac{d}{dt}\Bigg((1-t)^{\alpha+1}(1+t)^{\beta+1}\frac{d}{dt}y \Bigg)+\nonumber\\
 n(n+\alpha+\beta+1)(1-t)^{\alpha}(1-t)^{\beta}y=0\label{jacobsleq}
\end{eqnarray}
where $\alpha,\beta>-1$ (\cite{abr}). The solutions to this differential equation are denoted by $P^{(\alpha,\beta)}_n(t)$. Two solutions 
$P^{(\alpha,\beta)}_n$ are orthogonal over the interval $[-1,1]$ with weight factor $(1-t)^{\alpha}(1+t)^{\beta}$. 

The families of orthogonal polynomials 
such as Legendre and Chebyshev polynomials which are orthogonal over the interval $[-1,1]$ can be scaled appropriately so that one obtains corresponding 
families of orthogonal polynomials over the interval $[0,1]$. These are called shifted orthogonal polynomials (shifted Legendre or shifted Chebyshev).
 \section{Optimal steering using orthogonal polynomials}
 In this section, we show that for appropriately chosen cost functions on the input $\bold{u}$, orthogonal functions give optimal input to steer the nonholonomic 
 integrator from the origin to $(0,0,a)$ for some $a\in\mathbb{R}\setminus\{0\}$. 
 \begin{problem}\label{prob}
  We define the generalized quadratic cost on the inputs of the nonholonomic integrator
  \begin{equation}
 J=\frac{1}{2}\int_{t_0}^{t_f} \left[ \begin{array}{ccc} {u}_1(t)&{u}_2(t)\end{array}\right]
 \left[ \begin{array}{ccc} {a}_1(t)&0\\0&{a}_2(t)\end{array}\right]\left[ \begin{array}{ccc} {u}_1(t)\\{u}_2(t)\end{array}\right]dt\label{diagcost}
\end{equation}
${a}_1(t),{a}_2(t)>0$ $\forall t \in [t_0,t_f]$.
  The problem is to find the optimal input for $(\ref{noholintg})$ to steer the state from the origin to $(0,0,a)$ which minimizes $J$.
 \end{problem}
\begin{lemma}\label{lem1}
 Consider the nonholonomic integrator $(\ref{noholintg})$ with the following cost function on the input energy 
\begin{equation}
 J=\frac{1}{2}\int_{t_0}^{t_f} \left[ \begin{array}{ccc} {u}_1(t)&{u}_2(t)\end{array}\right]
 \left[ \begin{array}{ccc} {a}_1(t)&0\\0&{a}_2(t)\end{array}\right]\left[ \begin{array}{ccc} {u}_1(t)\\{u}_2(t)\end{array}\right]dt\nonumber
\end{equation}
where ${a}_1(t),{a}_2(t)>0$. Then the optimal input satisfies the following ode
\begin{equation}
 \dot{\Bigg(\left[ \begin{array}{c} a_1{u}_1\\a_2u_2\end{array} \right]\Bigg)}=\left[ \begin{array}{cc} 0&-2\lambda\\2\lambda&0\end{array} \right]
 \left[ \begin{array}{c} {u}_1\\u_2\end{array} \right] \label{optipdiagcost}
\end{equation}
where $\lambda\in\mathbb{R}$. 
\end{lemma}
\begin{IEEEproof}
 Using $\dot{x}_1=u_1$ and $\dot{x}_2=u_2$, the augmented Lagrangian is given by $L_a=1/2(a_1\dot{x}_1^2+a_2\dot{x}_2^2)+p(t)(\dot{x}_3-x_1\dot{x}_2+x_2\dot{x}_1)$. 
 Therefore,   
 \begin{eqnarray}
  &&\frac{\partial L_a}{\partial \dot{x}_1}=a_1\dot{x}_1+p(t)x_2,\; \frac{\partial L_a}{\partial {x}_1}=-p(t)\dot{x}_2\nonumber\\
  &&\frac{\partial L_a}{\partial \dot{x}_2}=a_2\dot{x}_2-p(t)x_1,\;\frac{\partial L_a}{\partial {x}_2}=p(t)\dot{x}_1\nonumber\\
  &&\frac{\partial L_a}{\partial \dot{x}_3}=p(t),\;\frac{\partial L_a}{\partial {x}_3}=0.\nonumber
 \end{eqnarray}
Applying Euler-Lagrange equations $(\ref{eullag})$, 
 \begin{eqnarray}
  \frac{d}{dt}({a_1\dot{x}_1}+p(t)x_2)&=&-p(t)\dot{x}_2\nonumber\\
  \frac{d}{dt}({a_2\dot{x}_2}-p(t)x_1)&=&p(t)\dot{x}_1\nonumber\\
  \frac{d}{dt}p(t)&=&0.\nonumber
 \end{eqnarray}
Let $p(t):=\lambda$ where $\lambda\in \mathbb{R}$. 
Now using $\dot{x}_1=u_1$ and $\dot{x}_2=u_2$, 
\begin{eqnarray}
 \frac{d}{dt}({a_1u_1})+2\lambda {u}_2=0\label{sl1}\\
 \frac{d}{dt}({a_2u_2})-2\lambda {u}_1=0\label{sl2}.
\end{eqnarray}
\end{IEEEproof}
\begin{remark}\label{sleq}
 From Equation $(\ref{sl1})$, $u_2=-\frac{1}{2\lambda}\frac{d}{dt}({a_1u_1})$ $\Rightarrow a_2u_2=-\frac{a_2}{2\lambda}\frac{d}{dt}({a_1u_1})$. Now 
substituting this in Equation $(\ref{sl2})$, 
\begin{equation}
 \frac{d}{dt}(-\frac{a_2}{2\lambda}\frac{d}{dt}({a_1u_1}))=2\lambda {u}_1 \Rightarrow 
 \frac{d}{dt}(a_2\frac{d}{dt}({a_1u_1}))=-4\lambda^2 {u}_1.\label{legpoly}
\end{equation}
This is a Sturm-Liouville equation whose solutions are orthogonal polynomials. Similarly, $u_2$ belongs to the family of solutions of 
\begin{equation}
 \frac{d}{dt}(a_1\frac{d}{dt}({a_2u_2}))=-4\lambda^2 {u}_2.\label{legpoly2}
\end{equation}

\end{remark}

 
Thus, both $u_1,u_2$ belong to a family of different orthogonal functions. But 
$u_2=-\frac{1}{2\lambda}\frac{d}{dt}({a_1u_1})$ may not be orthogonal to $u_1$ unless $a_1=a_2$. \\
{\bf Optimality of Chebyshev polynomials}: 
Observe that Chebyshev polynomials are orthogonal w.r.t. weight factor $\frac{1}{\sqrt{1-t^2}}$. 
Let $a_1(t)=a_2(t)={\sqrt{1-t^2}}$. Therefore, from $(\ref{legpoly})$,
\begin{equation}
 \frac{d}{dt}\Bigg(\sqrt{1-t^2}\frac{d}{dt}\Bigg((\sqrt{1-t^2})w_1\Bigg)\Bigg) =-4\lambda^2 {w}_1.\label{Chebyshevlike}
\end{equation}
Let $y=(\sqrt{1-t^2})w_1$. Therefore, 
\begin{eqnarray}
 \frac{d}{dt}\Bigg(\sqrt{1-t^2}\frac{d}{dt}y\Bigg) =-4\lambda^2 \frac{y}{\sqrt{1-t^2}}\nonumber\\
 \Rightarrow (1-t^2)\frac{d^2 y}{dt^2}-t\frac{dy}{dt}=-4\lambda^2 y.\label{Chebyshev}
\end{eqnarray}
The above ode is Chebyshev differential equation whose solutions are of the form 
\begin{equation}
 y=b_1T_{2\lambda}(t)+b_2\sqrt{1-t^2}U_{2\lambda-1}(t)\label{gensolcheb}
\end{equation}
where $T_{2\lambda}(t)$ is a Chebyshev polynomial of the first kind and $U_{2\lambda}(t)$ is a Chebyshev polynomial of the second kind (\cite{weis}). 
%
%
\begin{theorem}\label{chebopt}
 Consider the nonholonomic integrator $(\ref{noholintg})$ with the following cost function on the input energy 
 \begin{equation}
 J=\frac{1}{2}\int_{-1}^{1} \left[ \begin{array}{ccc} {u}_1&{u}_2\end{array}\right]
 \left[ \begin{array}{ccc} \sqrt{1-t^2}&0\\0&\sqrt{1-t^2}\end{array}\right]\left[ \begin{array}{ccc} {u}_1\\{u}_2\end{array}\right]dt.\label{diagcost1}
\end{equation}
Suppose that for a given state transfer, the optimal inputs exists and are continuously differentiable in time. 
Then they are given by solutions of the Sturm-Liouville equation $(\ref{Chebyshevlike})$.

\end{theorem}
\begin{IEEEproof}
 It follows by substituting $a_1=a_2=\sqrt{1-t^2}$ in Lemma $\ref{lem1}$ and Remark $\ref{sleq}$ that one obtains the Sturm-Liouville equation 
 $(\ref{Chebyshevlike})$ for $u_1$ and $u_2$. 
\end{IEEEproof}

To find the optimal inputs for a state transfer with the cost function given by $(\ref{diagcost1})$, 
we need to determine the Lagrange multiplier $\lambda$ which is an eigenvalue of the Sturm-Liouville equation $(\ref{gensolcheb})$. This is determined 
by the boundary conditions of the control problem. 

 Consider the steering problem (Problem $\ref{prob}$) for the nonholonomic integrator from $\bold{x}(-1) = (0,0,0)$ to 
 $\bold{x}(1) = (0,0,a)$, $a>0$ with 
 the cost function given by $(\ref{diagcost1})$.  
 We want to find the optimal inputs for this state transfer. It follows from Theorem $\ref{chebopt}$ that the optimal inputs are given by solutions 
 of Equation $(\ref{Chebyshevlike})$. In the following theorem we state how to compute these optimal solutions explicitly. 
 \begin{theorem}\label{chebcompute}
 Consider the steering problem for the nonholonomic integrator from $(0,0,0)$ at $t=-1$ to $(0,0,a)$ $(a>0)$ at $t=1$ with the 
 cost function given by $(\ref{diagcost1})$. Then, 
  the optimal inputs 
  are given by 
\begin{eqnarray}
u_1 = \frac{b_1T_{2\lambda}}{\sqrt{1-t^2}} + b_2U_{2\lambda - 1} \label{ip1}\\
u_2 = \frac{b_2T_{2\lambda}}{\sqrt{1-t^2}} - b_1U_{2\lambda - 1} \label{ip2}
\end{eqnarray}
where $\lambda=1$, $(b_1^2 + b_2^2)\frac{\pi}{2} = a$ and $T_{\lambda}$ and $U_{\lambda}$ are Chebyshev polynomials of the first and the second kind respectively. 
The optimal cost is equal to $a$. 
  \end{theorem}
\begin{IEEEproof}
 Refer Appendix $\ref{computn}$.
\end{IEEEproof}
The above theorem solves Problem $\ref{prob}$ for the cost function $(\ref{diagcost1})$. 

\begin{remark}
 Suppose we want the 
 state transfer along the negative $x_3$ direction i.e., $a<0$. Consider 
 the transformation $x_1'=x_2$, $x_2'=x_1$, $x_3'=x_3$ and $u_1'=u_2$, $u_2'=u_1$ in Lemma $\ref{chebcompute}$. Then, one obtains
 \begin{equation}
  \dot{x}_3'=x_2'u_1'-x_1'u_2'\nonumber
 \end{equation}
 with a boundary condition $x_3'(1)=a<0$ and $x_1'(1)=x_2'(1)=0$. Notice that it follows from the proof of Lemma $\ref{chebcompute}$ that 
 this transformation takes care of the minus sign appearing in the 
 $x_3$ coordinate and $u_1'=u_2$ and $u_2'=u_1$ are optimal inputs for a state transfer along the negative $x_3$ direction. 
\end{remark}
{\bf Optimality of Jacobi-like polynomials}: 
Note that Equations $(\ref{legpoly})$ and $(\ref{legpoly2})$ can be rewritten as
\begin{eqnarray}
 \frac{d}{dt}\{a_2\frac{d}{dt}k_1\}=-4\lambda^2 \frac{k_1}{a_1}\label{jac1}\\
 \frac{d}{dt}\{a_1\frac{d}{dt}k_2\}=-4\lambda^2 \frac{k_2}{a_2}\label{jac2}
\end{eqnarray}
where $k_1:=a_1u_1$ and $k_2:=a_2u_2$. Choosing 
\begin{eqnarray}
 &&a_1(t)=(1-t)^{-\alpha}(1-t)^{-\beta}, a_2(t)=(1-t)^{\alpha+1}(1+t)^{\beta+1},\nonumber\\
 &&4\lambda^2=n(n+\alpha+\beta+1)\label{jac3}
\end{eqnarray}
and comparing $(\ref{jac1})$ with $(\ref{jacobsleq})$, one obtains the Sturm-Liouville equation for the Jacobi polynomials. 
Therefore, 
\begin{equation}
 k_1=P_n^{(\alpha,\beta)}\label{jacsol}
\end{equation}
where $P_n^{(\alpha,\beta)}$ denotes the eigenfunction corresponding to $n$ for a fixed $\alpha,\beta$. 
Rewriting $(\ref{jac2})$ using $(\ref{jac3})$, 
\begin{eqnarray}
 \frac{d}{dt}\{(1-t)^{-\alpha}(1-t)^{-\beta}\frac{d}{dt}k_2\}=\nonumber\\
 -n(n+\alpha+\beta+1)(1-t)^{-\alpha-1}(1+t)^{-\beta-1} k_2.\nonumber
\end{eqnarray}
Let $-\alpha-1=\eta$ and $-\beta-1=\zeta$. Therefore, 
\begin{eqnarray}
 &&\frac{d}{dt}\{(1-t)^{\eta+1}(1-t)^{\zeta+1}\frac{d}{dt}k_2\}\nonumber\\
 &&=-n(n-\eta-\zeta-1)(1-t)^{\eta}(1+t)^{\zeta} k_2.\nonumber
\end{eqnarray}
Let $l=n-\eta-\zeta-1$. Therefore, 
\begin{eqnarray}
 \frac{d}{dt}\{(1-t)^{\eta+1}(1-t)^{\zeta+1}\frac{d}{dt}k_2\}=\nonumber\\
 -l(l+\eta+\zeta+1)(1-t)^{\eta}(1+t)^{\zeta} k_2 
 \Rightarrow k_2= P_l^{(\eta,\zeta)}.\label{jacsol1}
\end{eqnarray}
Thus, both $k_1$ and $k_2$ are Jacobi polynomials, hence, $u_1$ and $u_2$ are scaled versions of Jacobi polynomials. This leads to the following 
observation. 
\begin{theorem}
 Consider the nonholonomic integrator $(\ref{noholintg})$ with the cost function on the input energy given by $(\ref{diagcost})$. 
 Suppose $a_1(t),a_2(t)$ and $\lambda$ satisfy Equation $(\ref{jac3})$ where $-1<\alpha,\beta\le 0$. 
 Let $k_1=a_1u_1$ and $k_2=a_2u_2$ where $u_1,u_2$ are 
  the continuously differentiable optimal inputs for a given state transfer. 
  Then, $k_1,k_2$ are given by solutions of the Sturm-Liouville equation $(\ref{jacobsleq})$.
\end{theorem}
\begin{IEEEproof}
 It follows from $(\ref{jacobsleq})$ and $(\ref{jac1})$ that $k_1$ satisfies $(\ref{jacobsleq})$. 
 It also follows from the discussion above that if $-1<\alpha,\beta\le 0$, then $k_2$ satisfies $(\ref{jacobsleq})$.
\end{IEEEproof}
Notice that the optimal inputs are Jacobi polynomials scaled by appropriate factors i.e., $u_1=k_1/a_1$ and $u_2=k_2/a_2$ where $k_1,k_2$ are Jacobi polynomials. 
\begin{remark}
 When $\alpha=\beta=0$, one obtains Legendre polynomials. 
Therefore, from $(\ref{jac3})$, $a_1=1$ and $a_2=(1-t^2)$. Thus, $k_1$ is a Legendre polynomial and so is $u_1$. 
Notice that $\frac{d^2}{dt^2}k_2=-4\lambda^2k_2(1-t^2)^{-1}$. 
The solution set of this ode is not given by Legendre polynomials. 
Therefore, using $\alpha=\beta=0$ does not give optimal inputs where both the inputs are Legendre polynomials. 
For $\alpha=\beta=-\frac{1}{2}$, one obtains Chebyshev polynomials of the first kind whereas; for 
$\alpha=\beta=\frac{1}{2}$, one obtains Chebyshev polynomials of the second kind. 
\end{remark}
\begin{remark}
 As far as computations of the optimal input are concerned, one needs to find the eigenvalue $\lambda$ of the Sturm-Liouville operator using the 
 initial and terminal conditions of the state transfer. The corresponding eigenfunctions then give the optimal input. 

 Notice that although Legendre polynomials are not optimal for any of the cost function mentioned above, one can still do steering using only Legendre 
polynomials. Thus, apart from trigonometric polynomials, one can use Legendre polynomials, Chebyshev polynomials and Jacobi polynomials $(-1<\alpha, 
\beta\le0)$ to steer the nonholonomic integrator. 

\end{remark}

\section{Steering algorithm using orthogonal polynomials}
We now extend the steering algorithm using trigonometric polynomials given by Murray and Sastry (\cite{murray}) for families of orthogonal polynomials. 
However, \cite{murray} gave an algorithm for more general nonholonomic systems as well. In future, we want to consider steering of 
general nonholonomic systems using orthogonal polynomials. \\
{\bf Algorithm $1$}: 
\begin{itemize}
 \item Steer $x_1,x_2$ to their desired values using constant inputs. 
 \item Take a pair of orthogonal polynomials (trigonometric, Legendre, Chebyshev and so on) such that one of them is an even function and the other is an 
 odd function. Multiply by the weighting factor to guarantee the 
 orthogonality with the unit polynomial. Scale the polynomials appropriately to hit the desired point in the $x_3-$direction. The other two directions remain unchanged. 
\end{itemize}
\begin{remark}
 Note that if both $u_1$ and $u_2$ are even/odd functions as orthogonal polynomials over $[-1,1]$, then their respective integrals $x_1,x_2$ are 
 odd/even functions which makes $\dot{x}_3$ an odd function. Thus, $x_3(1)=0$. Therefore, while choosing orthogonal pairs, one must choose an even-odd pair 
 from any family of orthogonal polynomials. 
\end{remark}

\subsection{Steering using Legendre polynomials}
Notice that orthogonality of trigonometric polynomials was used in Example $\ref{sasex}$ to steer the state from the origin to $(0,0,a)$. We now demonstrate 
that one can exploit the orthogonality of Legendre polynomials to achieve the desired state transfer. 
\begin{example}
  Consider the shifted Legendre polynomials $1, 2t-1, \frac{1}{2}(3(2t-1)^2-1)$ which are orthogonal 
 on the interval $[0,1]$. Let $u_1=2t-1$ and 
 $u_2=6t^2-6t+1$. Thus, by orthogonality with $1$ on $[0,1]$, there is no motion in $x_1,x_2$ direction. Observe that 
 \begin{eqnarray}
  x_3&=&\int_0^1(t^2-t)(6t^2-6t+1)-(2t^3-3t^2+t)(2t-1)dt\nonumber\\
  &=&\int_0^1((6t^4-6t^3+t^2-6t^3+6t^2-t)-\nonumber\\
  &&(4t^4-6t^3+2t^2-2t^3+3t^2-t))dt\nonumber\\
  &=&\int_0^1 (2t^4-4t^3+2t^2)dt=\frac{2}{5}-1+\frac{2}{3}=\frac{1}{15}.\nonumber
 \end{eqnarray}
Thus, $u_1,u_2$ can be scaled so that $(0,0,a)$ can be reached. 
\end{example}
\subsection{Steering using Chebyshev polynomials of the first kind}
\begin{example}
 Consider the Chebyshev polynomials of the first kind $T_1(t)=t$ and $T_2(t)=2t^2-1$ which are orthogonal over $[-1,1]$. 
 Let $u_1(t)=t/\sqrt{1-t^2}$ and $u_2(t)=(2t^2-1)/\sqrt{1-t^2}$ which are normalized by the weighting factor for Chebyshev polynomials so that 
 $\int_{-1}^1u_i(t)=0$ since, $T_1,T_2$ are orthogonal to $1$ with respect to the weight factor $1/(\sqrt{1-t^2})$. 
 \begin{eqnarray}
  x_1(t)=\int_{-1}^t\frac{\tau}{\sqrt{1-\tau^2}}d\tau=\int_{-\frac{\pi}{2}}^{\sin ^{-1}t} \frac{\sin \theta}{\cos \theta}\cos \theta d\theta=\nonumber\\
  -[\cos (\sin ^{-1}t)-  \cos (\sin ^{-1}(-1))]=-\cos (\sin ^{-1}t)\nonumber\\
  x_2(t)=\int_{-1}^t\frac{2\tau^2-1}{\sqrt{1-\tau^2}}d\tau=-\int_{-\frac{\pi}{2}}^{\sin ^{-1}t} \frac{\cos 2\theta}{\cos \theta}\cos \theta d\theta=\nonumber\\
  -\frac{1}{2} [\sin (2\sin ^{-1}t)-  \sin (2\sin ^{-1}(-1))]=-\frac{1}{2} \sin (2\sin ^{-1}t).\nonumber
 \end{eqnarray}
Therefore, 
\begin{eqnarray}
 x_3(1)&=&\int_{-1}^1 [-\cos (\sin ^{-1}t)\frac{(2t^2-1)}{\sqrt{1-t^2}}+\nonumber\\
 &&\frac{1}{2} \sin (2\sin ^{-1}t)\frac{t}{\sqrt{1-t^2}}]dt\nonumber\\
 &=&\int_{-\frac{\pi}{2}}^{\frac{\pi}{2}}[\cos \theta\cos 2\theta+\frac{1}{2}\sin 2\theta\sin \theta]d\theta\nonumber\\
 &=&\frac{1}{2}\int_{-\frac{\pi}{2}}^{\frac{\pi}{2}}[\cos \theta\cos 2\theta+\cos\theta]d\theta\nonumber\\
 &=&\frac{1}{4}\int_{-\frac{\pi}{2}}^{\frac{\pi}{2}}[\cos \theta+\cos 3\theta + 2\cos \theta]d\theta\nonumber\\
 &=&\frac{1}{4}[3\sin \theta + \frac{1}{3}\sin \theta]|_{-\frac{\pi}{2}}^{\frac{\pi}{2}}=\frac{1}{2}[3+\frac{1}{3}]=\frac{5}{3}.\nonumber
\end{eqnarray}

 Thus, steering can be done in the $x_3-$direction.
 For steering along the $x_1-x_2$ plane, one can use constant polynomials which do not produce any motion along the $x_3-$direction. 
 One can similarly show that steering is possible with the Chebyshev polynomials of the second kind. 
\end{example}
\subsection{Simulations}
In this subsection, simulations are done using simulink of MATLAB of the nonholonomic integrator using two different kinds of orthogonal polynomials which are Legendre and Chebyshev(optimal steering) in time interval [-1,1] starting from (0,0,0) to (0,0,1).

For Legendre polynomials,  the inputs are chosen as 
$u_{1}(t) = \sqrt{\frac{15}{4}}P_{1}(t)$, $u_{2}(t) = \sqrt{\frac{15}{4}}P_{2}(t)$.
For Chebyshev polynomials, the inputs are chosen as
$u_{1}(t) = \sqrt{\frac{2}{\pi}}\frac{T_{2}(t)}{\sqrt{1-t^2}}$, $u_{2}(t) = -\sqrt{\frac{2}{\pi}}U_{1}(t)$.

The Results of Simulations are illustrated below (Note that Blue trajectories are due to Chebyshev inputs and Red trajectories are due to Legendre inputs).
%
%
%
\begin{figure}[ht]\label{plot1}
\begin{center}
\includegraphics[scale=0.35]{{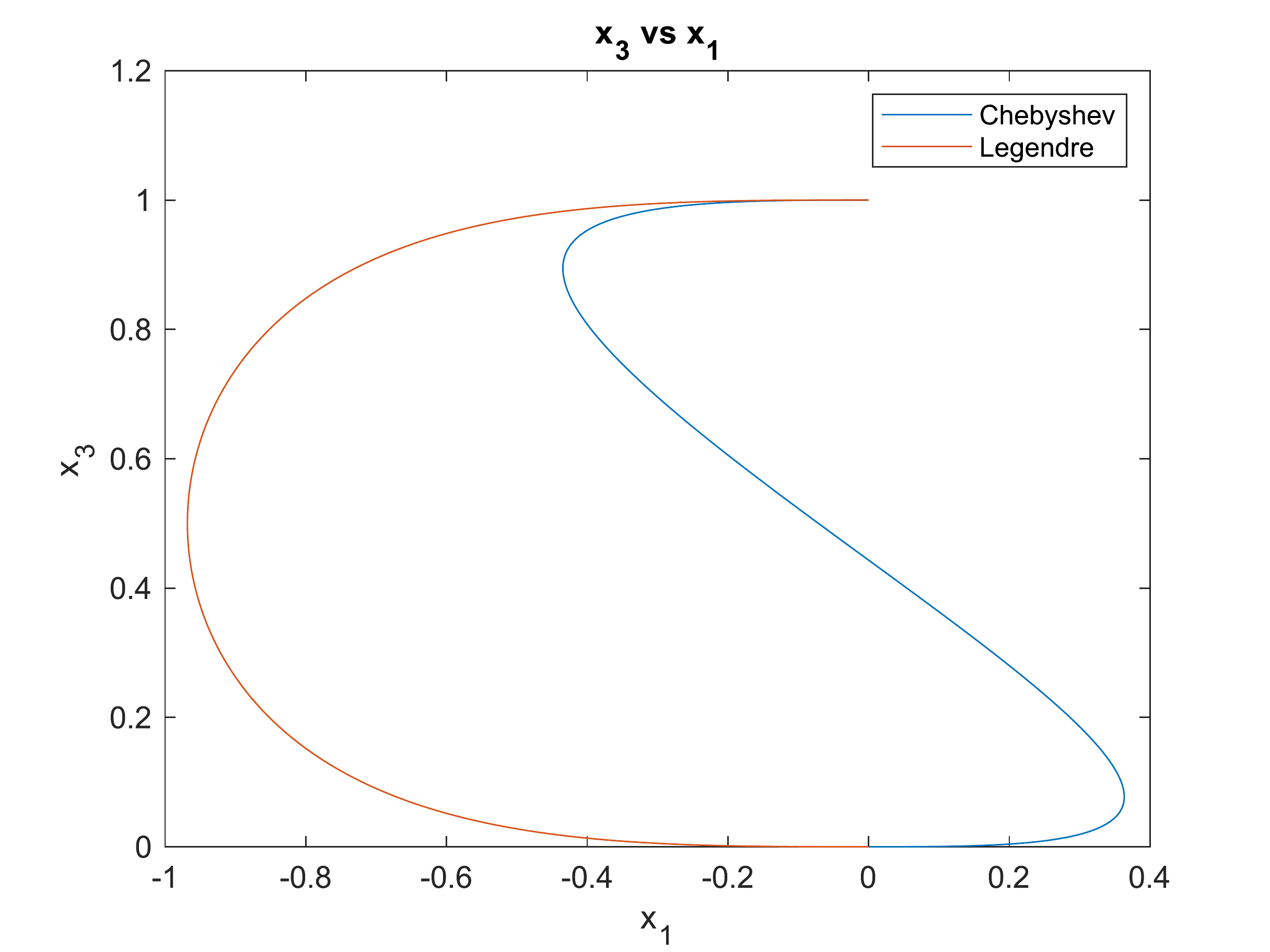}}
\caption{{Plot of $x_3$ vs $x_1$.}}
\end{center}
\end{figure} 
\begin{figure}[ht]\label{plot2}
\begin{center}
\includegraphics[scale=0.35]{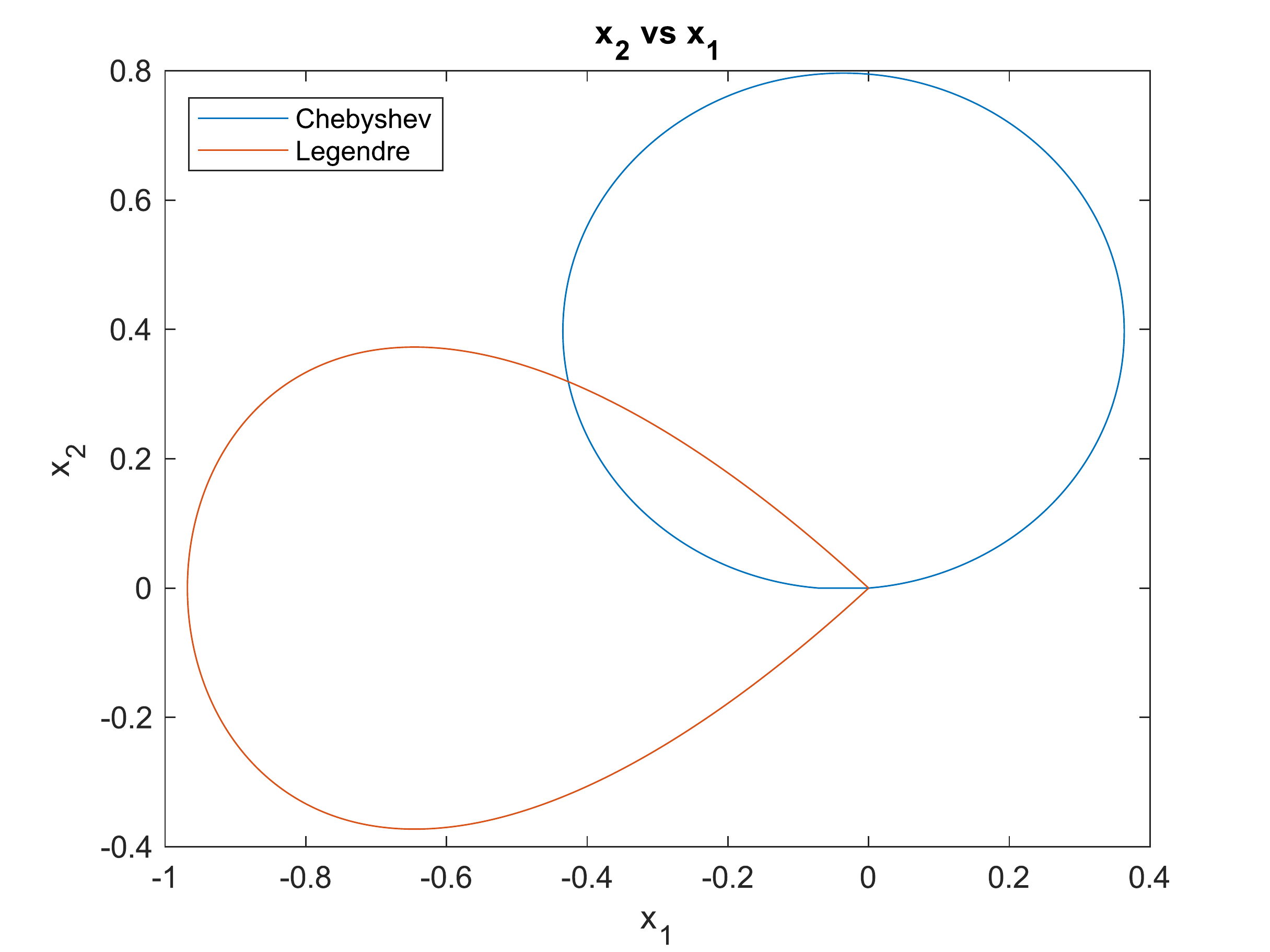}
\caption{{Plot of $x_2$ vs $x_1$.}}
\end{center}
\end{figure} 
\begin{figure}[ht]\label{plot4}
\begin{center}
\includegraphics[scale=0.35]{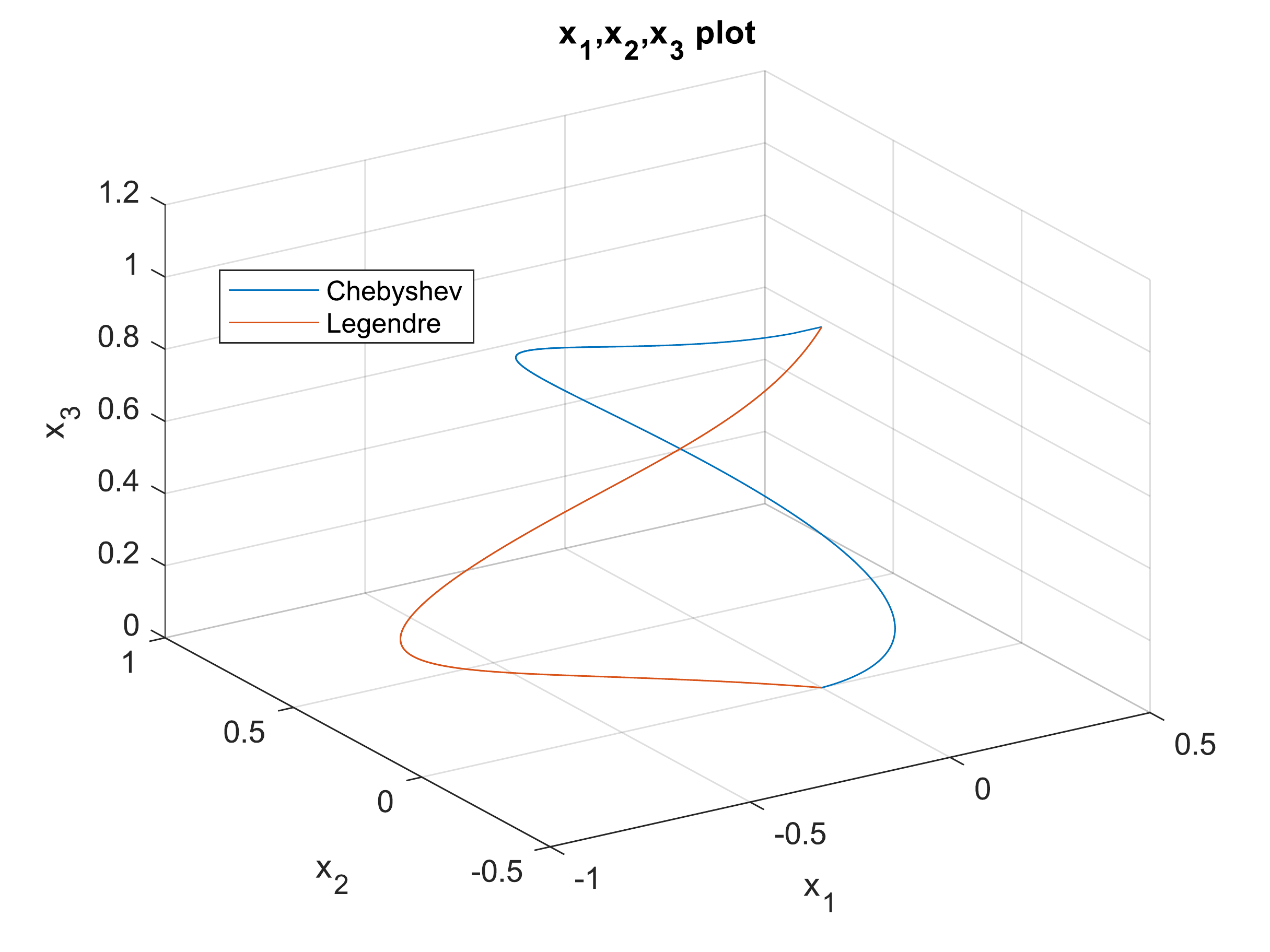}
\caption{{Evolution of the state trajectory.}}
\end{center}
\end{figure} 

\subsection{Steering the generalized nonholonomic integrator}
Consider the following generalization of $(\ref{noholintg})$
\begin{eqnarray}
 \dot{x}_i&=&u_i,\; i=1,\ldots,m,\nonumber\\
 \dot{x}_{ij}&=& x_iu_j-x_ju_i,\; i<j=1,\ldots,m\label{gennonhol} 
\end{eqnarray}
and consider the steering problem for this system. Murray and Sastry \cite{murray} gave a steering algorithm using sinusoids. We now show that one can do 
the steering using orthogonal polynomials as well. We exploit the fact that if $u_i,u_j$ $i\neq j$ are both even/odd orthogonal polynomials, then 
there is no steering in $x_i,x_j$ and $x_{ij}$ component. Steering in $x_{ij}$ happens when $u_i,u_j$ form an even-odd pair of orthogonal polynomials. \\
{\bf Algorithm $2$}: 
\begin{itemize}
 \item Choose $u_i,\; i=1,\ldots,m$ as constant polynomials and do the desired steering in $x_i, \; i=1,\ldots,m$. 
 \item For $i=1,\ldots,m-1$,\\
 \hspace*{.2in} For $j=1,\ldots,i-1$,\\
 \hspace*{.4in} $u_j\leftarrow 0$,\\
 \hspace*{.2in} End For\\
 \hspace*{.2in} $u_i \leftarrow$ any odd orthogonal polynomial (from any family). \\
 \hspace*{.2in} $u_k$ $\leftarrow$ any even orthogonal polynomial $(k=i+1\ldots,m)$.\\
 \hspace*{.2in} Steer  $x_{ik}$, $(k=i+1\ldots,m)$ to the desired value.\\
 End For.
\end{itemize}
The above algorithm ensures that $x_1,\ldots,x_m$ are steered first to the desired values, followed by $x_{12},\ldots,x_{1m}$, which is followed by 
$x_{23},\ldots,x_{2m}$ and so on up to 
the steering of $x_{m-1,m}$. This sequential steering is due to orthogonality properties of families of orthogonal polynomials and properties of 
even-odd functions where at each step, already steered variables remain undisturbed due to orthogonality and vanishing integrals of odd functions. 
\section{Control on $\mathbb{SO}(3)$ using orthogonal polynomials}
%
We now demonstrate the use of orthogonal polynomials for steering on the Lie group $\mathbb{SO}(3)$. 
We use Pontryagin's Maximum principle (PMP) on Lie groups and related results from \cite{Spindler} to obtain necessary conditions for optimal inputs. 
The Lie algebra for the Lie group $\mathbb{SO}(3)$ is given by $3\times 3$ skew symmetric matrices 
\begin{equation}
 {so}(3)=\Bigg\{\hat{\bold{\omega}}\in\mathbb{R}^{3\times3}\;|\;
 \hat{\bold{\omega}}=\left[ \begin{array}{ccc} 0&-\omega_3&\omega_2\\\omega_3&0&-\omega_1\\-\omega_2&\omega_1 &0\end{array}\right]
 \Bigg\} \nonumber
\end{equation}
which can be identified with $\mathbb{R}^3$ as $\bold{\omega}= \left[ \begin{array}{c} \omega_1\\\omega_2 \\\omega_3\end{array}\right]\mapsto 
\hat{\bold{\omega}}$. 

%
Consider the following control problem on the Lie group $\mathbb{SO}(3)$.
The attitude or the orientation of a spacecraft is a matrix $g\in \mathbb{SO}(3)$ whose columns form an orthonormal frame attached to the spacecraft. Let 
 $t\mapsto \omega(t)$ be the angular velocity of the space-craft. Let $\hat{\omega}$ be the skew symmetric matrix (defined above from a vector to get elements of Lie 
 algebra of $\mathbb{SO}(3)$) obtained from $\omega$. Then $\dot{g}=\hat{\omega}g$. Suppose 
 \begin{eqnarray}
  &&
  E_1=\left[ \begin{array}{cccc} 0&0&0\\0&0&-1\\0&1&0\end{array}\right], 
 E_2=\left[ \begin{array}{cccc} 0&0&-1\\0&0&0\\1&0&0\end{array}\right],\nonumber\\
 &&E_3=\left[ \begin{array}{cccc} 0&1&0\\-1&0&0\\0&0&0\end{array}\right].\nonumber
 \end{eqnarray}
This gives $\dot{g}=(\omega_1(t) E_1+\omega_2(t) E_2+\omega_3(t) E_3)g(t)$. 
A common problem in spacecraft attitude control is
now to perform an attitude maneuver which moves
the spacecraft from rest to rest between a given initial
attitude $g(t_0) = g_0$ and a prescribed target attitude
$g(t_1) = g_1$. Since it is sensible to try to perform such
a maneuver while keeping the overall angular velocities
low, we consider a cost functional of the form $\int_{t_0}^{t_1}(\omega_1^2(t)+\omega_2^2(t)+\omega_3^2(t))dt$ which we want to minimize. 
We refer the reader to \cite{Spindler} and the references therein where this problem was studied before. The problem is 
\begin{eqnarray}
 &\mbox{minimize}& \int_{t_0}^{t_1}(\omega_1^2(t)+\omega_2^2(t)+\omega_3^2(t))dt\nonumber\\
 &\mbox{subject to}& \dot{g}=Ug, g(t_0)=g_{0}, g(t_1)=g_1\label{liegr} 
\end{eqnarray}
where $U=\sum_{i=1}^{3}\omega_iE_i$. Let ${p}$ denote co-states 
which acts on the elements of the Lie algebra as follows 
\begin{eqnarray}
 {p}(E_i)=p_i,\;i=1,2,3.\nonumber
\end{eqnarray}
Using PMP on Lie groups (\cite{Spindler}), it turns out that the Lagrange multipliers ${p}$ satisfy  
\begin{equation}
 \dot{{p}}=-{p}\circ \mbox{ad}(U).\label{costeq}
\end{equation}
 Moreover, 
\begin{equation}
 [E_1,E_2]=E_3,\;[E_2,E_3]=E_1,\;[E_3,E_1]=E_2.\label{liebrak}
\end{equation}
Therefore, using $(\ref{costeq})$ and $(\ref{liebrak})$, one obtains 
\begin{eqnarray}
 \dot{p}_1&=&\dot{p}(E_1)=-p \circ([U,E_1])=-p \circ(-\omega_2E_3+\omega_3E_2)\nonumber\\
 &=&\omega_2p_3-\omega_3p_2\label{spattctrl1}\\
 \dot{p}_2&=&\dot{p}(E_2)=-p \circ([U,E_2])=-p \circ(\omega_1E_3-\omega_3E_1)\nonumber\\
 &=&\omega_3p_1-\omega_1p_3\label{spattctrl2}\\
 \dot{p}_3&=&\dot{p}(E_3)=-p \circ([U,E_3])=-p \circ(-\omega_1E_2+\omega_2E_1)\nonumber\\
 &=&\omega_1p_2-\omega_2p_1.\label{spattctrl3}
\end{eqnarray}

The Hamiltonian is $H= 
p_1\omega_1+p_2\omega_2+p_3\omega_3-(\omega_1^2+\omega_2^2+\omega_3^2)$ (\cite{Spindler}). 
Using the first order necessary conditions for maximization of $H$, $\frac{\partial{H}}{\partial \omega_i}=0$ $\Rightarrow$
$2\omega_i=p_i$. Substituting in $(\ref{spattctrl1})-(\ref{spattctrl3})$, we obtain 
\begin{eqnarray}
 2\dot{\omega}_1&=& 2\omega_2\omega_3-2\omega_3\omega_2=0,\nonumber\\ 2\dot{\omega}_2&=& -2\omega_1\omega_3+2\omega_3\omega_1=0,\nonumber\\
 2\dot{\omega}_3&=& 2\omega_1\omega_2-
 2\omega_2\omega_1=0.\nonumber
\end{eqnarray}
Therefore, $\omega_i=c_i$, where the constants $c_i$ can be obtained from the boundary conditions. 

If the cost function is $\int_{t_0}^{t_1}q(t)(\omega_1^2(t)+\omega_2^2(t)+\omega_3^2(t))dt$, then one obtains $\omega_i={c_i}/{q}$ 
(\cite{Spindler}). Suppose $\omega_3=1$. 
Therefore, $\dot{g}=(\omega_1(t) E_1+\omega_2(t) E_2+ E_3)g(t)$. Let the cost function be $\frac{1}{2}\int_{t_0}^{t_1}(\omega_1^2(t)+\omega_2^2(t))dt$. 
Using $\dot{{p}}=-{p}\circ$ad$(U)$ and Equations $(\ref{spattctrl1})-(\ref{spattctrl3})$, 
\begin{equation}
 \dot{p}_1= \omega_2p_3-p_2,\;\dot{p}_2= -\omega_1p_3+p_1,\;\dot{p}_3= \omega_1p_2-\omega_2p_1.\label{spattctrl4}
\end{equation}
From the maximization property of $H=p_1\omega_1+p_2\omega_2+p_3-(\omega_1^2+\omega_2^2)$, using the first order conditions $\frac{\partial{H}}{\partial \omega_i}=0$, one obtains $p_1=2\omega_1$, 
$p_2=2\omega_2$ and $p_3=2c$. Therefore, from $(\ref{spattctrl4})$, 
\begin{eqnarray}
 \dot{\omega}_1=\dot{p}_1=\omega_2(c-1), \; \dot{\omega}_2=\dot{p}_2=-\omega_1(c-1)\nonumber\\
 \Rightarrow \ddot{\omega}_1=-(c-1)^{2}\omega_1, \; \ddot{\omega}_2=-(c-1)^{2}\omega_2.\label{optctrlso3}
\end{eqnarray}
Thus, the optimal inputs are given by sinusoids. 

Now suppose the cost function is $\frac{1}{2}\int_{t_0}^{t_1}q(t)(\omega_1^2(t)+\omega_2^2(t))dt$. Proceeding as done in the previous paragraph, one obtains 
$p_1=q\omega_1$, $p_2=q\omega_2$ and $p_3=c$. Therefore, from $(\ref{spattctrl4})$, we obtain 
\begin{eqnarray}
 \frac{d}{dt}(q\omega_1)=\dot{p}_1=\omega_2(c-q)\nonumber\\
 \frac{d}{dt}(q\omega_2)=\dot{p}_2=-\omega_1(c-q)\nonumber\\
 \Rightarrow \frac{d}{dt}\left[ \begin{array}{cccc} q\omega_1\\q\omega_2\end{array}\right]=
 \left[ \begin{array}{cccc} 0&c-q\\-(c-q)&0\end{array}\right]\left[ \begin{array}{cccc} \omega_1\\\omega_2\end{array}\right]. \label{optctrlso3wc1}
\end{eqnarray}
Therefore, it follows that 
\begin{eqnarray}
 \frac{d}{dt}\Bigg(\frac{q}{(c-q)}\frac{d}{dt}(q\omega_1)\Bigg)=\frac{d}{dt}(q\omega_2)=-(c-q)\omega_1\label{slso31}\\
 \frac{d}{dt}\Bigg(\frac{q}{(c-q)}\frac{d}{dt}(q\omega_2)\Bigg)=-\frac{d}{dt}(q\omega_1)=-(c-q)\omega_2.\label{slso32}
\end{eqnarray}
Thus, one obtains a Sturm-Liouville equation whose solutions give optimal inputs. These inputs functions are orthogonal polynomials.
\begin{remark}
 Optimal control on the special Unitary group $\mathbb{SU}(2)$ is considered in \cite{Spindler} as an example of control of a quantum spin system. One can use the above approach to obtain a 
 Sturm-Liouville equation for an optimal control problem on $\mathbb{SU}(2)$ with two inputs having the cost function of the type 
 $\frac{1}{2}\int_{t_0}^{t_1}q(t)(u_1^2(t)+u_2^2(t))dt$. 
\end{remark}

\section{Sub-optimal fuel minimization using orthogonal polynomials}
%
Consider the fuel minimization problem 
\begin{equation}
 J=\int_{t_0}^{t_1}(|u_1(t)|+|u_2(t)|)dt \label{fuelmin}
\end{equation}
for the nonholonomic integrator. Suppose one wants to steer $(\ref{noholintg})$ from the origin to $(0,0,a)$. Any even-odd pair of orthogonal polynomials (trigonometric, 
Legendre, Chebyshev, Jacobi) work. Moreover, one can choose a set of orthogonal polynomials say Legendre polynomials and any even-odd pair of Legendre polynomials 
$(P_i(t),P_j(t))$ where $P_i(t),P_j(t)\neq 1$, $i\neq j$ can do the state transfer when scaled appropriately. Let $u_1(t)=b_1P_1(t)$ and $u_2=b_2P_2(t)$ where 
$P_1(t)$ and $P_2(t)$ denote the first and the second Legendre polynomial over $[-1,1]$. Then, one can choose 
$b_1,b_2$ to minimize $(\ref{fuelmin})$. 
Thus, we have the following optimization problem 
\begin{eqnarray}
 &\mbox{minimize} & \int_{t_0}^{t_1}(|b_1||P_1(t)|+|b_2||P_2(t)|)dt\nonumber\\
 &\mbox{subject to}& a= \int_{t_0}^{t_1}(x_2(t)u_1(t)-x_1(t)u_2(t))dt.\label{optprob1}
\end{eqnarray}
Let $c_i=\int_{t_0}^{t_1}P_i(t)dt$, $i=1,2$. Using $x_j(t)=b_j\int_{t_0}^tP_{j}(\tau)d\tau$ $j=1,2$, one can rewrite $(\ref{optprob1})$ as 
\begin{eqnarray}
 &\mbox{minimize} & |b_1|c_1+|b_2|c_2\nonumber\\
 &\mbox{subject to}& c= b_1b_2\label{optprob2}
\end{eqnarray}
where $c$ is an appropriate constant. 
Now by the AM-GM inequality on positive real numbers, we get  
\begin{equation}
    |b_1|c_1+|b_2|c_2 \geq 2\sqrt{|c|c_1c_2}.
\end{equation}
As the inequality is sharp, the minimum of $J$ is $2\sqrt{|c|c_1c_2}$ i.e.,
\begin{equation}\label{sopfuel}
    \mbox{min}(J) = 2\sqrt{|c|c_1c_2}.
\end{equation}
\begin{example}
Consider $u_1(t)$ = $b_1 P_1(t)$ and $u_2(t)$ = $b_2 P_2(t)$, where $P_n(t)$ denote the $n-$th Legendre polynomial. 
Also consider different inputs as $u_1(t) = d_1\sin(\pi t)$ and $u_2(t) = d_2\cos(\pi t)$. 
We compare the values of $J$ for these inputs. From Equation (\ref{sopfuel}), we just need to compute $|c|$, $c_1$ and $c_2$ for the given sets of inputs, 
in the interval $[-1,1]$ for state transfer from $(0,0,0)$ to $(0,0,1)$.

For  $u_1(t)$ = $b_1 P_1(t)$, $u_2(t)$ = $b_2 P_2(t)$, we obtain  
\begin{eqnarray}
 c_1 = \int_{-1}^{1}|P_1(t)|dt,
 c_2 = \int_{-1}^{1}|P_2(t)|dt, \nonumber\\
 c = \frac{-1}{2\int_{-1}^{1}P_2(t)\big(\int_{-1}^{\tau}P_1(t)dt\big) d\tau}.\nonumber
\end{eqnarray}
Computing these integrals gives us 
\begin{eqnarray}
c_1 = 1, c_2 = 0.7698, c = 3.75 \nonumber\\
\mbox{min}(J) = 3.3981.\nonumber
\end{eqnarray}
For  $u_1(t)$ = $b_1 \sin(\pi t)$, $u_2(t)$ = $b_2 \cos(\pi t)$, we obtain  
\begin{eqnarray}
 c_1 = \int_{-1}^{1}|\sin(\pi t)|dt, 
 c_2 = \int_{-1}^{1}|\cos(\pi t)|dt, \nonumber\\
 c = \frac{-1}{2\int_{-1}^{1}\cos(\pi t)\big(\int_{-1}^{\tau}\sin(\pi t)dt\big) d\tau}.\nonumber
\end{eqnarray}
Computing these integrals gives us 
\begin{eqnarray}
c_1 = 1.2732, c_2 = 1.2732, c = 3.1407\nonumber\\
\mbox{min}(J) = 4.5135.\nonumber
\end{eqnarray}
Thus, in this example, Legendre polynomials are better than trigonometric functions at optimizing the cost function given by the $\mathcal{L}_1$ norm 
of the input.
\end{example}
\begin{remark}
 One can use the same approach to find sub-optimal inputs from a family of orthogonal functions when $\mathcal{L}_p$ norm of the input 
is used as a cost function.
\end{remark}


\section{Conclusion}
We showed that for the nonholonomic integrator, families of orthogonal polynomials such as Chebyshev, Legendre and Jacobi polynomials can do the 
steering in addition to the trigonometric polynomials. Furthermore, we showed that for an appropriately defined cost 
function on the inputs, the optimal inputs are orthogonal polynomials which are solutions of an appropriate Sturm-Liouville differential equation. 
We showed that for some specific state transfers and a specific cost function, Chebyshev polynomials give optimal inputs. 
We showed that the steering algorithm of \cite{murray} can be extended 
for families of orthogonal polynomials for the nonholonomic integrator and the generalized nonholonomic integrator. We also showed that for an under-actuated optimal control on the Lie group 
$\mathbb{SO}(3)$, optimal inputs are given by a certain Sturm-Liouville equation. Furthermore, we showed that using orthogonal polynomials, one 
can construct sub-optimal solutions for each family of orthogonal polynomials by solving a finite dimensional optimization problem. 

In future, we want to extend these ideas to the optimal steering of generalized nonholonomic integrator, extended nonholonomic integrator and general 
nonholonomic systems. We also want to investigate is it possible to construct a cost function for which Legendre polynomials give optimal solutions. 
\appendices
\section{Properties of Chebyshev polynomials}\label{chebprop}
Let $T_n(t)$ and $U_n(t)$ denote Chebyshev polynomials of first and second kind respectively. We list all properties 
(\cite{mason}, Chapter $1$ and $2$) which are used in the proof of Theorem $\ref{chebcompute}$ in Appendix $\ref{computn}$.
\begin{enumerate}
	\item\label{Ch1} Let $(1-t^2) {d^2 y \over d t^2} - x {d y \over d t} + p^2 y = 0$ be the Chebyshev differential equation. 
	Then, its general solution for $t\in[-1,1]$ is given by $$ y = b_1 T_p(t) + b_2\sqrt{1-t^2}U_{p-1}(t).   $$
	\item \label{Ch2} $T_n(cos\theta) = cos(n\theta)$, $U_n(cos\theta) = \frac{sin\big((n+1)\theta\big)}{sin\theta}$.
	\item \label{Ch3} $T_n(-t) = (-1)^n T_n(t)$, $U_n(-t) = (-1)^n U_n(t)$. That is, Chebyshev polynomials of even order are even functions 
	and Chebyshev polynomials of odd order are odd functions.
	\item \label{Ch4} $\frac{d T_n}{d t} = n U_{n - 1}$. $\frac{d U_n}{d t} = \frac{(n + 1)T_{n + 1} - t U_n}{t^2 - 1}$.
	\item \label{Ch5} $\int U_n\ dt = \frac{T_{n + 1}}{n + 1}$, $\int T_n\, dt = \frac{1}{2} \left(\frac{T_{n + 1}}{n + 1} - \frac{T_{n - 1}}{n - 1}\right)$.
	\item \label{Ch6} $\int_{-1}^1 T_n(t)T_m(t)\,\frac{dt}{\sqrt{1-t^2}}=
	\begin{cases}
	0 & n\ne m \\
	\pi & n=m=0 \\
	\frac{\pi}{2} & n=m\ne 0.
	\end{cases}$
	\item \label{Ch7} $\int_{-1}^1 U_n(x)U_m(t)\sqrt{1-t^2}\,dt =
	\begin{cases}
	0 & \text{if } n\ne m, \\
	\frac{\pi}{2} & \text{if } n=m.
	\end{cases}$
\end{enumerate}
\section{Proof of Theorem \ref{chebcompute}}\label{computn}
It follows from Theorem $\ref{chebopt}$ that the optimal inputs are given by solutions of Equation $(\ref{Chebyshevlike})$. 
It follows from Equations $(\ref{legpoly})$ and $(\ref{legpoly2})$, both $u_1$ and $u_2$ are solutions of the same Sturm-Liouville equation since 
$a_1=a_2=\sqrt{1-t^2}$. 
From Equations 
$(\ref{gensolcheb})$ and $(\ref{Chebyshevlike})$, it follows that one of the optimal inputs say $u_1$ must be of the form $(\ref{ip1})$. 
Using Equation $(\ref{sl1})$ and $(\ref{ip1})$, 
\begin{eqnarray}
 2\lambda u_2=-\frac{d}{dt}(au_1)=-\frac{d}{dt}\Bigg({b_1T_{2\lambda}} + b_2{\sqrt{1-t^2}}U_{2\lambda - 1}\Bigg).\nonumber
\end{eqnarray}
Using Property $(\ref{Ch4})$ of Appendix $\ref{chebprop}$, 
\begin{eqnarray}
 2\lambda u_2&=&-(2\lambda b_1U_{2\lambda-1}-b_2\frac{t}{\sqrt{1-t^2}}U_{2\lambda-1}+\nonumber\\
 &&b_2\sqrt{1-t^2}\frac{(2\lambda T_{2\lambda}-tU_{2\lambda})}{t^2-1})\nonumber\\
 &=&-2\lambda b_1U_{2\lambda-1}+2\lambda b_2\frac{T_{2\lambda}}{\sqrt{1-t^2}}\nonumber
\end{eqnarray}
which implies that $u_2$ satisfies $(\ref{ip2})$. Therefore, 
\begin{eqnarray}
u_1 = \frac{b_1T_{2\lambda}}{\sqrt{1-t^2}} + b_2U_{2\lambda - 1} \nonumber\\
u_2 = \frac{b_2T_{2\lambda}}{\sqrt{1-t^2}} - b_1U_{2\lambda - 1}. \nonumber
\end{eqnarray}


From the terminal conditions $x_1(1)=x_2(1)=0$, it follows that
\begin{eqnarray}
\int_{-1}^1 u_1 d\tau = \int_{-1}^1 u_2 d\tau = 0. \label{ipcdn}
\end{eqnarray}
Substituting expressions for $u_1$ and $u_2$ 
in $(\ref{ipcdn})$ and using the orthogonality property (Property $(\ref{Ch6})$ of Chebyshev polynomials 
of the first kind from Appendix $\ref{chebprop}$), we obtain 
 $\int_{-1}^1b_2U_{2\lambda - 1}d\tau=0$ and $\int_{-1}^1b_1U_{2\lambda - 1}d\tau=0$. Now, using Property $(\ref{Ch5})$ of Chebyshev polynomials of the 
 second kind from Appendix $\ref{chebprop}$, 
 \begin{eqnarray}
  \int_{-1}^1U_{2\lambda - 1}d\tau=\frac{T_{2\lambda}(1)-T_{2\lambda}(-1)}{2\lambda}.\label{chebcdn}
 \end{eqnarray}
Using Property $(\ref{Ch2})$ from Appendix $\ref{chebprop}$, $T_{2\lambda}(1)=T_{2\lambda}(\cos 0)=\cos(2\lambda. 0)=1$ and 
similarly, $T_{2\lambda}(-1)=T_{2\lambda}(\cos \pi)=\cos( 2\lambda\pi)$. Now to satisfy $(\ref{ipcdn})$, 
it follows from $(\ref{chebcdn})$ that $T_{2\lambda}(1)-T_{2\lambda}(-1)=0$ which implies that $\cos( 2\lambda\pi)=1$. Therefore, $\lambda$
must be an integer. Note that we have used the terminal conditions $x_1(1)=0$ and $x_2(1)=0$ to conclude this. 

We now want to find $\lambda$ from the given terminal conditions on $x_3$ i.e., $x_3(1)=a$. 
Notice that 
\begin{eqnarray}
\dot{x}_3 = x_1\dot{x}_2 - x_2\dot{x}_1 \nonumber\\
\dot{x}_3 = \frac{d}{dt}(x_1x_2) - 2x_2\dot{x}_1 \nonumber\\
x_3(1) = -2\int_{-1}^1 x_2u_1 d\tau = a.
\end{eqnarray}
Now, by substituting expressions of $u_1$ and $x_2$ we get,
\begin{eqnarray}
-2\int_{-1}^1 u_1(\tau)(\int_{-1}^\tau u_2(t) dt) d\tau = a \Rightarrow\nonumber\\
\int_{-1}^1 (\frac{b_1T_{2\lambda}}{\sqrt{1-\tau^2}} + b_2U_{2\lambda - 1})(\int_{-1}^\tau \frac{b_2T_{2\lambda}}{\sqrt{1-t^2}} - b_1U_{2\lambda - 1} dt) d\tau = -\frac{a}{2}.\label{pt0} 
\end{eqnarray}
Notice that in the above integral, 
\begin{eqnarray}
\int_{-1}^\tau (\frac{b_2T_{2\lambda}}{\sqrt{1-t^2}} - b_1U_{2\lambda - 1}) dt =  \int_{-1}^\tau \frac{b_2T_{2\lambda}}{\sqrt{1-t^2}} dt \nonumber\\
-\frac{b_1}{2\lambda}[T_{2\lambda}(\tau) -T_{2\lambda}(-1)] \label{pt1}
\end{eqnarray}
where we have used Property $(\ref{Ch5})$ of Chebyshev polynomials of the second kind from Appendix $\ref{chebprop}$. Note that $\lambda \neq0$ from Equation $(\ref{pt1})$. 
Now
\begin{eqnarray}
\int_{-1}^\tau \frac{b_2T_{2\lambda}}{\sqrt{1-t^2}} dt = \int_{\cos^{-1}\tau}^\pi \frac{b_2T_{2\lambda}(\cos\theta)}{\sin\theta} \sin\theta d\theta \nonumber\\
= b_2\int_{\cos^{-1}\tau}^\pi \cos(2\lambda\theta) d\theta = -\frac{b_2}{2\lambda}\sin(2\lambda\cos^{-1}\tau).\label{pt2} 
\end{eqnarray}
where we have used Property $(\ref{Ch2})$ of Chebyshev polynomials from Appendix $\ref{chebprop}$. Substituting $(\ref{pt2})$ in $(\ref{pt1})$, 
\begin{eqnarray}
 \int_{-1}^\tau (\frac{b_2T_{2\lambda}}{\sqrt{1-t^2}} - b_1U_{2\lambda - 1}) dt =  -\frac{b_2}{2\lambda}\sin(2\lambda\cos^{-1}\tau)\nonumber\\
 -\frac{b_1}{2\lambda}[ T_{2\lambda}(\tau) - T_{2\lambda}(-1)].\label{pt3}
\end{eqnarray}
Substituting $(\ref{pt3})$ in $(\ref{pt0})$, 
\begin{eqnarray}
 &&\int_{-1}^1 (\frac{b_1T_{2\lambda}(\tau)}{\sqrt{1-\tau^2}} + b_2U_{2\lambda - 1}(\tau))(\int_{-1}^\tau \frac{b_2T_{2\lambda}}{\sqrt{1-t^2}} - 
 b_1U_{2\lambda - 1} dt) d\tau
 = \nonumber\\
&& \int_{-1}^1 (\frac{b_1T_{2\lambda}(\tau)}{\sqrt{1-\tau^2}} + b_2U_{2\lambda - 1}(\tau))( -\frac{b_1}{2\lambda}[ T_{2\lambda}(\tau) - T_{2\lambda}(-1)]\nonumber\\
&& -\frac{b_2}{2\lambda}\sin(2\lambda\cos^{-1}\tau))d\tau=\nonumber\\
&&   \int_{-1}^1(\frac{-b_1^2T_{2\lambda}^2(\tau)}{2\lambda\sqrt{1-\tau^2}} + \frac{b_1^2T_{2\lambda}(\tau)T_{2\lambda}(-1)}{2\lambda\sqrt{1-\tau^2}} 
-\nonumber\\
&&  \frac{b_1b_2T_{2\lambda}(\tau)}{2\lambda\sqrt{1-\tau^2}}\sin(2\lambda\cos^{-1}\tau))  -\frac{b_1b_2}{2\lambda}U_{2\lambda - 1}(\tau)T_{2\lambda}(\tau) + \nonumber\\
   &&\frac{b_1b_2}{2\lambda}U_{2\lambda - 1}(\tau)T_{2\lambda}(-1) -\frac{b_2^2}{2\lambda}U_{2\lambda - 1}(\tau)\sin(2\lambda\cos^{-1}\tau) )d\tau.\label{pt6}
\end{eqnarray}
Note that using Property $(\ref{Ch3})$ of Chebyshev polynomials from Appendix $\ref{chebprop}$, the third, the fourth and the fifth term in the above integral form odd functions, hence, 
vanish. Moreover, using the orthogonality property (Property $(\ref{Ch6})$ of Chebyshev polynomials from Appendix $\ref{chebprop}$) of $T_{2\lambda}(\tau)$ with $T_0(\tau)=1$, 
the integral of the second term also vanishes. Therefore, $(\ref{pt6})$ can be simplified as, 
\begin{eqnarray}
 &&\int_{-1}^1 (\frac{b_1T_{2\lambda}(\tau)}{\sqrt{1-\tau^2}} + b_2U_{2\lambda - 1}(\tau))(\int_{-1}^\tau \frac{b_2T_{2\lambda}}{\sqrt{1-t^2}} - 
 b_1U_{2\lambda - 1} dt) d\tau
 = \nonumber\\
 &&\int_{-1}^1(\frac{-b_1^2T_{2\lambda}^2(\tau)}{2\lambda\sqrt{1-\tau^2}} -\frac{b_2^2}{2\lambda}U_{2\lambda - 1}(\tau)\sin(2\lambda\cos^{-1}\tau) )d\tau.\label{pt7}
\end{eqnarray}
The first term in  the above integral can be simplified using Property $(\ref{Ch2})$ of Chebyshev polynomials of the first kind 
(Appendix $\ref{chebprop}$) as
\begin{eqnarray}
 &&\int_{-1}^1\frac{-b_1^2T_{2\lambda}^2(\tau)}{2\lambda\sqrt{1-\tau^2}}d\tau=
 \int_{0}^{\pi}\frac{-b_1^2\cos^2(2\lambda\theta)}{2\lambda\sin \theta}\sin \theta d\theta\nonumber\\
 &&=-\frac{\pi b_1^2}{4\lambda}.\label{pt9}
\end{eqnarray}
Notice that the second term in $(\ref{pt7})$ can be simplified as
\begin{eqnarray}
 &&\int_{-1}^1\frac{b_2^2}{2\lambda}U_{2\lambda - 1}(\tau)\sin(2\lambda\cos^{-1}\tau) d\tau=\nonumber\\
 &&\int_{0}^\pi \frac{b_2^2}{2\lambda}\big(U_{2\lambda - 1}(\cos\theta)\big)\sin(2\lambda\theta) \sin\theta d\theta=\nonumber\\
 &&\int_{0}^\pi \frac{b_2^2}{2\lambda} \sin^2(2\lambda\theta) d\theta = - \frac{\pi b_2^2}{4\lambda}\label{pt8}
\end{eqnarray}
where we used Property $(\ref{Ch2})$ of Chebyshev polynomials of the second kind $U_{2\lambda - 1}(\cos\theta)=\sin(2\lambda\theta)/\sin(\theta)$. 
Substituting $(\ref{pt9})$ and $(\ref{pt8})$ in $(\ref{pt7})$,
\begin{eqnarray}
 &&\int_{-1}^1 (\frac{b_1T_{2\lambda}(\tau)}{\sqrt{1-\tau^2}} + b_2U_{2\lambda - 1}(\tau))(\int_{-1}^\tau \frac{b_2T_{2\lambda}}{\sqrt{1-t^2}} - 
 b_1U_{2\lambda - 1} dt) d\tau\nonumber\\
 &=&-\frac{\pi(b_1^2 + b_2^2)}{4\lambda}.\label{pt10}
\end{eqnarray}
Comparing $(\ref{pt0})$ and $(\ref{pt10})$, 
\begin{equation}
 \frac{\pi(b_1^2 + b_2^2)}{4\lambda}=\frac{a}{2}\Rightarrow (b_1^2 + b_2^2)\frac{\pi}{2} = \lambda a.\label{pt11}
\end{equation}

Now consider the cost function $J$, 
by substitution of inputs $u_1,u_2$ and by using the properties of Chebyshev polynomials, 
\begin{eqnarray}
J &=& \int_{-1}^1 {\sqrt{1-t^2}}(u_1^2 + u_2^2) dt \\
u_1^2 &=& \bigg(\frac{b_1T_{2\lambda}}{\sqrt{1-t^2}} + b_2U_{2\lambda - 1}\bigg)^2 \nonumber\\
&=& \frac{b_1^2T_{2\lambda}^2}{{1-t^2}} + b_2^2U_{2\lambda - 1}^2 + 2\frac{b_1b_2T_{2\lambda}}{\sqrt{1-t^2}}U_{2\lambda - 1}  \\
u_2^2 &=& \bigg(\frac{b_2T_{2\lambda}}{\sqrt{1-t^2}} - b_1U_{2\lambda - 1}\bigg)^2 \nonumber\\
&=& \frac{b_2^2T_{2\lambda}^2}{{1-t^2}} + b_1^2U_{2\lambda - 1}^2 - 2\frac{b_1b_2T_{2\lambda}}{\sqrt{1-t^2}}U_{2\lambda - 1}  \\
\Rightarrow J &=& \int_{-1}^1 (b_1^2 + b_2^2)(\frac{T_{2\lambda}^2}{\sqrt{1-t^2}} + U_{2\lambda - 1}^2\sqrt{1-t^2}) dt \label{pt12}\\ 
 &=& (b_1^2 + b_2^2)\frac{\pi}{2} = \lambda a\label{pt13}
\end{eqnarray}
where the last equation follows from using properties $(\ref{Ch6})$ and $(\ref{Ch7})$ of Chebyshev polynomials from Appendix $\ref{chebprop}$ 
and Equation $(\ref{pt11})$. 
Since $J\ge 0$, $\lambda\ge0$.  
Thus, $\lambda = 1$ for minimizing $J$ since $\lambda\neq0$ and from Equation $(\ref{pt13})$, $(b_1^2 + b_2^2)\frac{\pi}{2} = a$. 
The optimal inputs are given by Equations $(\ref{ip1})$ and $(\ref{ip2})$ for $\lambda=1$ 
and the optimal cost is $J=a$. This completes the proof.
\section*{Acknowledgement}
Authors are thankful to K.V. Manohar for help in simulations.
\bibliographystyle{ieeetr}        
\bibliography{steering}


\end{document}